\documentclass[12pt]{amsart}
\usepackage{amscd,amsthm,amsfonts,amssymb,amsmath,euscript}
\usepackage[dvips]{graphicx}
\usepackage[dvips]{graphics}
\usepackage[T2A]{fontenc}
\usepackage[cp866]{inputenc}

\newtheorem{theorem}{Theorem}[section]
\newtheorem{lemma}{Lemma}[section]

\newtheorem{note}{Remark}[section]

\newcommand{\Aut}{\mathop{\sf Aut}\nolimits}

\begin{document}

\centerline {\bf Brane Topological Field Theories }

\centerline {\bf and Hurwitz numbers for CW-complexes}

{\ }

\centerline {\bf Sergey M. Natanzon{\footnote{Supported by grants
RFBR-08-01-00110, NWO 047.011.2004.026 (05-02-89000-NWO-a),
NSh-709.2008.1}}}

{\ }

\centerline{A.N.Belozersky Institute, Moscow State University}

\centerline{Independent University of Moscow}

\centerline{Institute Theoretical and Experimental Physics}

\centerline{natanzon@mccme.ru}

{\ }

\smallskip

\section{Introduction}\label{r0}

We expand Topological Field Theory of \cite{At},\cite{W} on some
special CW-complexes (brane complexes). This Brane Topological Field
Theory extends the Open Topological Field Theory from
\cite{L},\cite{MS},\cite{AN}. For arbitrary two-dimension brane
complexes the Brane Topological Field Theory is the "open part" of
Cyclic Foam Topological Field Theory from \cite{N}. We prove that
the Brane Topological Field Theory one-to-one corresponds to
infinite dimensional  Frobenius Algebras, graduated by CW-complexes
of lesser dimension.

We define general and regular Hurwitz numbers of brane complexes and
prove that they generate Brane Topological Field Theories. For
general Hurwitz numbers corresponding algebra is an algebra of
coverings of lesser dimension. For regular Hurwitz numbers the
Frobenius algebra is an algebra of families of subgroups of finite
groups.

\smallskip

\section{Brane Topological Field Theory}\label{r1}

\subsection{Simple CW-complexes} \label{r1.1}
A finite compact CW-complex is said to be \textit{simple} if
\begin{itemize}
\item it is the closure of theirs cells of maximal dimension;
\item the closure of any cell is a closed ball.
\end{itemize}
\textit{An isomorphism} between simple CW-complexes is a
homeomorphism conserving theirs cell structures. Zero-dimensional
cells are called \textit{vertexes}. Denote by $\Omega_b$ the set of
all vertexes of $\Omega$.

Denote by $|M|$ the cardinality of a set $M$. Following by Smith and
Dold (see also \cite{BR}) we define \textit{$d$-brunched covering}
as a continuous map $h: X\rightarrow Y$ between two Hausdorff spaces
and a continuous map $t: Y\rightarrow \textmd{Sym}^d(X)$ such that:
1) $x\in th(x)$ for all $x\in X$, 2) $\textmd{Sym}^d(h)(ty)=dy$ for
all $y\in Y$. We say that a $d$-brunched covering is \textit{finite}
if $Y$ assume a CW-complex structure such that the function
$|h^{-1}(y)|$ is constant on any cell.

\begin{lemma}\label{l1.2}  Any simple CW-complexes has
a finite brunched covering over a ball.
\end{lemma}

\textsl{Proof}. Let $\Omega$ be a simple CW-complex with vertexes
$\{q_k|k\in K\}$. Consider a ball $B$ and $k$ different points
$\{b_k\in\partial B|k\in K\}$. Then there exists a finite number
isomorphic classes of simple CW-complexes with vertexes $\{b_k|k\in
K\}$ that form a stratification of $B$. Consider a set $\{B_i|i\in
I\}$ of representatives of theirs isomorphic classes.

Let $\{\Omega_j|j\in J\}$ be the set of closing of
$\dim\Omega$-dimensional cells $\omega_j\subset\Omega$. Then for any
$j\in J$ there exists an isomorphism of simple CW-complexes
$h_j:\Omega_j\rightarrow B_{i(j)}$. Gluing now the manifolds
$\{\Omega_j|j\in J\}$ along its boundary by maps $h_k^{-1}h_j$,
using the same rules that manifolds $\{\Omega_j|j\in J\}$ form
$\Omega$. Then we get a simple CW-complex $\widetilde{\Omega}$ that
isomorphic to $\Omega$. The maps $h_j$ generate a finite branching
covering $h:\widetilde{\Omega}\rightarrow B$.

$\Box$
\begin{theorem}\label{t1.1}A Hausdorff space $\Omega$ has a
stratification, generating a simple CW-complex if and only if there
exists its finite branching covering over a ball
$h:\Omega\rightarrow B$.
\end{theorem}

\textsl{Proof}. Let $h:\Omega\rightarrow B$ be a finite branching
covering. Construct a stratification of $B$ onto cells, such that,
the functions $|h^{-1}(y)|$ is constant on any cell. Preimages of
these cells form a structure of simple CW-complex on $\Omega$. The
inverse affirmation follows from lemma \ref{l1.2}.

$\Box$

\subsection{Brane complexes} \label{r1.2}
A simple CW-complex is said to be \textit{oriented} if all cells of
positive dimensional are oriented. Fix a set $S$ of \textit{colors}.
An oriented complex $\Omega$ is called \textit{colored complex} if a
color $s(l)\in S$ is assigned to any its cell $l$ of positive
dimension, and for any connected component of $\Omega$ any two
celles of the same dimension have different colors. \textit{An
isomorphism} between colored complexes is an isomorphism between
simple CW-complexes, conserving the orientations and the colors of
cells.

A system of cyclic orders on vertexes of connected components of a
oriented complex $\Omega$ is called a \textit{cyclic order on}
$\Omega$. A colored complex with a cyclic order is called
\textit{cyclic complex}.

A connected subcomplex $\gamma\subset\Omega$ of codimension 1 into
cyclic complex $\Omega$ is called a \textit{cut} if:
\begin{itemize}
\item the restriction of $\gamma$ to any closed cell of $\omega\in\Omega$
either is empty, or homeomorphic to closed ball of codimension 1
in $\omega$;
\item $\gamma$ divides a connected component of $\Omega$ into two
connected components that splits its vertices into two nonempty
groups, compatible with the cyclic order on $\Omega$.
\end{itemize}

A cyclic complex is called a \textit{brane complex}, if for any
compatible with the cyclic order division of vertexes $\Omega$ there
exists a cut that realize it.

The structure of CW-complex on $\Omega$ forms a structure of
CW-complex on any its cut. A cut $\sigma_q$, dividing a vertex $q$
from other vertexes, is called \textit{vertex complex}. It separates
a neighborhood $U_q$ of $q$, that is a cone over $\sigma_q$. Any
cell $\sigma\in\sigma_q$ divides an oriented cell $\omega\in\Omega$.
Equip $\sigma$ by orientation, as a boundary of connected component
of $\omega\setminus\sigma$ that don't contain $q$. The colore
structure on $\Omega$ generate a colore structure on $\sigma_q$.
Thus, a brane complex generates a structure of connected colored
complexes on its vertex complexes.

\subsection{Topological Field Theory}\label{r1.3}
Below we assume that all vector spaces are defined over a field
$\mathbb{K}\supset\mathbb{Q}$. Let $\{X_m | m\in M\}$ be a finite
set of $n=|M|$ vector spaces $X_m$ over the field of complex numbers
$\mathbb{C}$. The action of the symmetric group $S_n$ on
$\{1,\dots,n\}$ induces its action on the sum of the vector spaces
$\left(\oplus_{\sigma} X_{\sigma(1)}\otimes\dots\otimes
X_{\sigma(n)}\right)$, where $\sigma$ runs over the bijections
$\{1,\dots, n\}\to M$, an element $s\in S_n$ takes
$X_{\sigma(1)}\otimes\dots\otimes X_{\sigma(n)}$ to
$X_{\sigma(s(1))}\otimes\dots\otimes X_{\sigma(s(n))}$. Denote by
$\otimes_{m\in M} X_m$ the subspace of all invariants of this
action. The vector space $\otimes_{m\in M} X_m$ is canonically
isomorphic to the tensor product of all $X_m$ in any fixed order;
the isomorphism is the projection of $\otimes_{m\in M} X_m$ to the
summand that is equal to the tensor product of $X_m$ in that order.

Denote by $\Sigma=\Sigma(D,S)$ the set of all isomorphism classes of
D-dimensional connected colored complexes. The inversion of the
orientations generates the involution $*:\Sigma\rightarrow\Sigma$.
Denote it by $\sigma\mapsto\sigma^*$.

Consider a family of finite-dimensional vector spaces
$\{B_{\sigma}|\sigma\in\Sigma\}$ and a family of tensors
$\{K_{\sigma}^\otimes\in B_{\sigma}\otimes
B_{\sigma^*}|\sigma\in\Sigma\}$. Using these data, we define now a
functor $\mathcal{V}$ from the monoidal category $\mathcal{S}$ of
$(D+1)$-dimensional brane complexes to the category of vector
spaces. This functor assigns the vector space $V_\Omega=
(\otimes_{q\in\Omega_b}B_q)$ to any $(D+1)$-dimensional brane
complex $\Omega$. Here $B_q$ is the copy of $B_{\sigma_q}$ that is a
vector space with a fixed isomorphism $B_q\rightarrow B_{\sigma_q}$.

\smallskip

We are going to describe all morphisms of the monoidal category
$\mathcal{S}$ and morphisms of the category of vector spaces that
correspond to it. \vskip0.3cm (1) {\it Isomorphism.} Define an
isomorphism of brane complexes as an isomorphism $\phi:\Omega\to
\Omega'$ of colored complexes, preserving the cyclic orders. Define
$\mathcal{V}(\phi)=\phi_*: V_{\Omega}\to V_{\Omega'}$ as linear
operator generated by the bijections
$\phi|_{\Omega_b}:\Omega_b\to\Omega'_b$. \vskip0.3cm
(2)\textit{Cut.} Let $\Omega$ be a connected brane complex and
$\gamma\subset\Omega$ be a cut. The cut $\gamma$ is represented by
two complexes $\gamma_+$ and $\gamma_-$ on the closure
$\overline{\Omega\setminus\gamma}$ of $\Omega\setminus\gamma$ .
Contract these complexes to points $q_+=q_+[\gamma]$ and
$q_-=q_-[\gamma]$, respectively. The contraction produces a brane
complex $\Omega'=\Omega[\gamma]$. Its vertices
$\Omega'=\Omega[\gamma]$ are the vertices of $\Omega$ and the points
$q_+$, $q_-$. The cyclic order, orientation and the coloring of
$\Omega$ induce an orientation and a coloring of $\Omega'$. Thus, we
can assume that $\Omega'$ is a brane complex and
$V_{\Omega'}=V_{\Omega}\otimes B_{q_+}\otimes B_{q_-}$. The functor
takes the morphism $\mathcal{V}(\eta)(x)=\eta_*(x)=x\otimes
K_{\sigma}^\otimes$, where $\sigma=\sigma_{q_+}=\sigma_{q_-}^*$, to
the morphism $\eta:\Omega\to\Omega'$. \vskip0.3cm(3) The tensor
product in $\mathcal{S}$ defined by the disjoint union of brane
complexes $\Omega'\otimes\Omega''\to \Omega'\coprod\Omega''$ induces
the tensor product of vector spaces $\theta_*: V_{\Omega'}\otimes
V_{\Omega''}\to V_{\Omega'\sqcup\Omega''}$.

The functorial properties of $\mathcal{V}$ can be easily verified.

\smallskip

Fix a tuple of vector spaces and vectors
$\{B_{\sigma},K_{\sigma}^\otimes\in B_{\sigma}\otimes
B_{\sigma^*}|\sigma\in\Sigma\}$, defining the functor $\mathcal{V}$.
A family of linear forms $\mathcal {F}= \{\Phi_\Omega:V_{\Omega}\to
\mathbb{K}\}$ defined for all brane complexes $\Omega\in\mathcal{S}$
is called a \textit{Brane Topological Field Theory} if it satisfies
the following axioms:

\vskip0.6cm $1^\circ$ {\it  Topological invariance.}
$$\Phi_{\Omega'}(\phi_*(x))=\Phi_\Omega(x)$$ for any isomorphism
$\phi:\Omega\to \Omega'$ of brane complexes.

\vskip 0.6cm $2^\circ$ {\it  Non-degeneracy.}\vskip 0.3cm Let
$\Omega$ be a brane complex with only two vertices $q_1$, $q_2$.
Then $\sigma_{q_2}=\sigma^*$ if $\sigma_{q_1}=\sigma$. Denote by
$(.,.)_\sigma$ the bilinear form $(.,.)_\sigma: B_\sigma\times
B_{\sigma^*}\rightarrow\mathbb{K}$, where
$(x',x'')_\sigma=\Phi_\Omega(x'_{q_1}\otimes x''_{q_2})$. Axiom
$2^\circ$ asserts that the forms $(.,.)_\sigma$ are non-degenerated
for all $\sigma\in\Sigma$.

\vskip 0.6cm $3^\circ$ {\it Cut invariance.}
$$\Phi_{\Omega'}(\eta_*(x))=\Phi_\Omega(x)$$
for any cut morphism $\eta:\Omega\to \Omega'$ of brane complexes.
\vskip 0.6cm

$4^\circ$ {\it Multiplicativity.}

$$\Phi_{\Omega}(\theta_*(x'\otimes x'))= \Phi_{\Omega'}(x')\Phi_{\Omega''}(x'')$$

for $\Omega=\Omega'\coprod\Omega''$,  $x'\in V_{\Omega'}$, $x''\in
V_{\Omega''}$.

\smallskip

Note that a Topological Field Theory defines the tensors
$\{K_{\sigma}^\otimes\in B_{\sigma}\otimes
B_{\sigma^*}|\sigma\in\Sigma\}$, since it is not difficult to prove:
\begin{lemma} \label{l1}Let $\{\Phi_\Omega\}$ be a Brane Topological Field Theory.
Then $(K_\sigma^\otimes,x_1 \otimes x_2)_\sigma= (x_1,x_2)_\sigma$,
for all $x_1\in B_\sigma$, $x_2\in B_{\sigma^*}$.
\end{lemma}

\begin{note} \label{note} Brane complexes with single cell of higher dimensional
form subcategory into the subcategory (ball category) of brane
complexes. Thus, we can consider a restriction Brane Topological
Field Theory on ball category.
\end{note}

\section{Colored Frobenius Algebras}\label{r2}

\subsection{Algebra} \label{r2.1}

We say that a connected brane complex $\Omega$ is a
\textit{compatible complex} for connected colored complexes
$\sigma_1,\sigma_2,...,\sigma_n$ if $\{\sigma_i\}$ are vertex of
$\Omega$ and the numeration of $\{\sigma_i\}$ generates the cyclic
order on $\Omega$. Denote by
$\Omega(\sigma_1,\sigma_2,...,\sigma_n)$ the set of all isomorphism
classes of compatible complexes for
$\sigma_1,\sigma_2,...,\sigma_n$. Then
$\Omega(\sigma_1,\sigma_2,...,\sigma_n)$ is either empty or consists
of a single element.

Let $\Omega(\sigma_1,\sigma_2,\sigma_3,\sigma_4)\neq\emptyset$. Then
there exist unique classes of cuts
$\sigma_{(1,2|3,4)},\sigma_{(4,1|2,3)}\in\Sigma$ such that
$\Omega(\sigma_1,\sigma_2,\sigma_{(1,2|3,4)})\neq\emptyset$,
$\Omega(\sigma_{(3,4|1,2)},\sigma_3,\sigma_4)\neq\emptyset$,
$\Omega(\sigma_4,\sigma_1,\sigma_{(4,1|2,3)})\neq\emptyset$,
$\Omega(\sigma_{(4,1|2,3)},\sigma_2,\sigma_3)\neq\emptyset$ and
$\sigma_{(3,4|1,2)}=\sigma_{(1,2|3,4)}^*$,
$\sigma_{(2,3|4,1)}=\sigma_{(4,1|2,3)}^*$.

{\ }

Consider a tuple of finite dimensional vector spaces
$\{B_{\sigma}|\sigma\in\Sigma\}$. Its direct sum $B_\star =
\bigoplus_{\sigma\in\Sigma}B_{\sigma}$ is called a \textit{colored
vector space}. A colored vector space with a bilinear form
$(.,.):B_\star\times B_\star\rightarrow\mathbb{K}$ and a tree-linear
form $(.,.,.):B_\star\times B_\star\times
B_\star\rightarrow\mathbb{K}$ is called a \textit{ Colored Frobenius
Algebra} if
\begin{itemize}
\item the form $(.,.)$  is not-degenerate;
\item $(B_{\sigma_1},B_{\sigma_2}) = 0$ for
$\sigma_1\neq\sigma_2^*$;
\item $(B_{\sigma_1},B_{\sigma_2},B_{\sigma_3}) = 0$ for
$\Omega(\sigma_1,\sigma_2,\sigma_3)=\emptyset$
\item $\sum_{i,j}(x_1,x_2,b_i^{(1,2|3,4)})F^{ij}_{(1,2|3,4)}(b_j^{(3,4|1,2)},x_3,x_4)=
\sum_{i,j}(x_4,x_1,b_i^{(4,1|2,3)})F^{ij}_{(4,1|2,3)}(b_j^{(2,3|4,1)},x_2,x_3)$.
\end{itemize}
Here $x_k\in B_{\sigma_k}$, $\{b_i^{(s,t|k,r)}\}$ is a basis of
$B_{\sigma_{(s,t|k,r)}}$ and $F^{ij}_{(s,t|k,r)}$ is the inverse
matrix for $(b_i^{(s,t|k,r)},b_j^{(k,r|s,t)})$.

{\ }

We will consider $B_*$ as an algebra with the multiplication
$(x_1x_2,x_3) =(x_1,x_2,x_3)$, for $x_k\in B_{\sigma_k}$. The axiom
$\sum_{i,j}(x_1,x_2,b_i^{(1,2|3,4)})F^{ij}_{(1,2|3,4)}(b_j^{(3,4|1,2)},x_3,x_4)=
\sum_{i,j}(x_4,x_1,b_i^{(4,1|2,3)})F^{ij}_{(4,1|2,3)}(b_j^{(2,3|4,1)},x_2,x_3)$
is equivalent to associativity for the algebra $B_*$. Moreover
$B_*$ is a Frobenius Algebra in the sense of \cite{F} if its
dimension is finite.

\subsection{Topology-algebra correspondence}  \label{r2.2}
\begin{theorem} \label{t2.1} Let $\mathcal {F}=
\{\Phi_\Omega:V_{\Omega}\to \mathbb{K}\}$ be a Brane Topological
Field Theory on a tuple of finite-dimensional vector spaces
$\{B_{\sigma}|\sigma\in\Sigma\}$. Then the poli-linear forms
\begin{itemize}
\item $(x',x'')=\Phi_{\Omega(\sigma_1,\sigma_2)}(x'_{q_1}\otimes
x''_{q_2})$, where $x'\in B_{\sigma_1}$, $x''\in B_{\sigma_2}$
\item $(x',x'',x''')=\Phi_{\Omega(\sigma_1,\sigma_2,\sigma_3)}
(x'_{q_1}\otimes x''_{q_2}\otimes x'''_{q_3})$, где $x'\in
B_{\sigma_1}$, $x''\in B_{\sigma_2}$, $x'''\in B_{\sigma_3}$.
\end{itemize}
generate a structure of Colored Frobenius Algebra on $B_\star =
\bigoplus_{\sigma\in\Sigma} B_{\sigma}$.
\end{theorem}

\textsl{Proof}. Only the last axiom is not obvious. Let us consider
a brane complex
$\Omega\in\Omega(\sigma_1,\sigma_2,\sigma_3,\sigma_4)$, and a cut
between the pairs of vertexes $\sigma_1,\sigma_2$ and
$\sigma_3,\sigma_4$. Then the cut-invariant axiom  and lemma
\ref{l1} give
$\sum_{i,j}(x_1,x_2,b_i^{(1,2|3,4)})F^{ij}_{(1,2|3,4)}(b_j^{(3,4|1,2)},x_3,x_4)
= \Phi_\Omega (x_1,x_2,x_3,x_4)$. Similarly,
$\sum_{i,j}(x_4,x_1,b_i^{(4,1|2,3)})F^{ij}_{(4,1|2,3)}(b_j^{(2,3|4,1)},x_2,x_3)
= \Phi_\Omega (x_1,x_2,x_3,x_4)$

$\Box$

\begin{theorem} \label{t2.2} Let $B_\star = \bigoplus_{\sigma\in\Sigma}
B_{\sigma}$ be a Colored Frobenius Algebra with poli-linear forms
$(.,.)$ and $(.,.,.)$. Then it generates a Brane Topological Field
Theory on $\{B_{\sigma}|\sigma\in\Sigma\}$ by means of following
construction. Fix a basis $\{b_i^\sigma\}$ of any vector space
$B_\sigma$, $\sigma\in\Sigma$. Consider the matrix $F^{ij}_\sigma$
that is the inverse matrix for
$F_{ij}^\sigma=(b_i^\sigma,b_j^{\sigma^*})$. Define the linear
functionals for any connected brane complexes by
\begin{itemize}
\item $\Phi_{\Omega(\sigma_1,\sigma_2...,\sigma_n)}(x^1_{q_1}\otimes
x^2_{q_2}\otimes ...\otimes
x^n_{q_n})=\sum_{\varsigma_1,\varsigma_2,...,\varsigma_{n-3}\in\Sigma}
(x^1,x^2,b_{i_1}^{\varsigma_1})F_{i_1j_1}^{\varsigma_1}$
$(b_{j_1}^{\varsigma_1^*},x^3_{q_3},b_{i_2}^{\varsigma_2})\\F_{i_2j_2}^{\varsigma_2}
(b_{j_2}^{\varsigma_2^*},x^4_{q_4},b_{i_3}^{\varsigma_3})...
...F_{i_{n-4}j_{n-4}}^{\varsigma_{n-4}}
(b_{j_{n-4}}^{\varsigma_{n-4}^*},x^{n-2}_{q_{n-2}},b_{i_{n-3}}^{\varsigma_{n-3}})
F_{i_{n-3}j_{n-3}}^{\varsigma_{n-3}}
(b_{j_{n-3}}^{\varsigma_{n-3}^*},x^{n-1}_{q_{n-1}},x^n_{q_n})$,\\
where $x^i_{q_i}\in B_{\sigma_i}.$
\end{itemize}
Then the family $\{\Phi_{\Omega(\sigma_1,\sigma_2...,\sigma_n)}\}$
generates a Brane Topological Field Theory.
\end{theorem}

\textsl{Proof}. The  topological invariance follows from the
invariance under cyclic renumbering the vertices of $\Omega$. The
invariance under the renumbering $q_i\mapsto q_j$  $j\equiv
i+1(\mathrm{mod} 2)$ follows from the last axiom for the tree-linear
form. The cut invariance follows directly from the definition of
$\Phi$ if we renumber the vertices marking the cut divide the
vertices $q_1,q_2,...q_k$ and $q_{k+1},q_{k+2},...q_n$.

$\Box$

These two theorems determine the one-to-one correspondence between
Brane Topological Field Theories and isomorphic classes of Colored
Frobenius Algebras.

\section{Hurwitz numbers of colored complexes}

\subsection{Hurwitz numbers}
Classical Hurwitz numbers are weighted number of brunched coverings
over closed surface with prescribed types of critical values
\cite{H}. Let us define some analog of  Hurwitz numbers for brane
complexes.

A brunched covering between simple CW-complexes
$f:\widetilde{\Omega}\rightarrow\Omega$ is called \textit{simple
covering of degree $d$} or \textit{simple $d$-covering} it is a
local homeomorphism on any cell, preimage of any cell consists of
unions of cells, and the number of cells from preimage is $d$ for
any cell of higher dimension. Simple coverings
$f':\widetilde{\Omega}'\rightarrow\Omega$ and
$f":\widetilde{\Omega}"\rightarrow\Omega$ are called
\textit{isomorphic} if there exists an isomorphism
$\widetilde{\varphi}:\widetilde{\Omega}'\rightarrow\widetilde{\Omega}"$
such that $f'=f"\widetilde{\varphi}$. Denote by $\Aut(f)$ the
automorphisms group of a simple covering $f$.

Let $\Omega'$ and $\Omega"$ colored complexes. Simple coverings
$f':\widetilde{\Omega}'\rightarrow\Omega'$ and
$f":\widetilde{\Omega}"\rightarrow\Omega"$ are called
\textit{equivalent} if there exists isomorphisms
$\widetilde{\varphi}:\widetilde{\Omega}'\rightarrow\widetilde{\Omega}"$
and $\varphi:\Omega'\rightarrow\Omega"$ such that $\varphi
f'=f"\widetilde{\varphi}$. Denote by $\Upsilon_{\Omega}^d$ the set
of equivalent classes of simple $d$-coverings over a colored complex
$\Omega$.

A simple $d$-covering $f:\widetilde{\Omega}\rightarrow\Omega$ over a
brane complex $\Omega$ generate a simple $d$-covering
$f_q:\widetilde{\Omega}_q\rightarrow\sigma_q$ over the vertex
complex $\sigma_q$ of a vertex $q$. The equivalence class of $f_q$
is called a \textit{local invariant} of $f$ at $q$. Consider a set
$\{\beta_q\}=\{\beta_q\in\Upsilon_{\sigma_q}^d|q\in\Omega_b\}$.
Denote by $\Upsilon_{\Omega}^d(\{\beta_q\})$ the set of isomorphic
classes of simple $d$-coverings
$f:\widetilde{\Omega}\rightarrow\Omega$ with local invariants
$f_q=\beta_q$ for all vertexes $q\in\Omega_b$. The weighted number
of simple $d$-covering
$$H^d(\Omega,\{\beta_q\}) = \sum_{f\in
\Upsilon_{\Omega}^d(\{\beta_q\})} 1/|\Aut(f)|$$ is called
\textit{Hurwitz number of simple $d$-coverings}

\subsection{Hurwitz Brane Topological Field Theory}
For any $\sigma\in\Sigma=\Sigma(D,S)$, denote by $C_{\sigma}^d$ the
vector space generated by the set $\Upsilon_{\sigma}^d$. The
involution $\sigma\mapsto\sigma^*$ generates the involution
$C_{\sigma}^d\mapsto C_{\sigma^*}^d$. Hurwitz numbers generate a
family of linear functions $\mathcal {H}^d=
\{\Phi_\Omega:V_{\Omega}\to \mathbb{K}\}$ on
$\{V_{\Omega}=(\otimes_{q\in\Omega_b}C_q^d)\}$, depending from brane
complexes $\Omega$.

\begin{theorem} \label{t3.1}  $\mathcal
{H}^d$ is a Brane Topological Field Theory.
\end{theorem}
\textsl{Proof}. It is not evident only the cut-invariant axiom. Fix
a cut $\sigma\subset\Omega\in\Omega(\sigma_1,...,\sigma_m)$ and a
covering $\beta\in\Upsilon_{\sigma}^d$. Let $\Omega\setminus\sigma=
\Omega_+\cup\Omega_-$. Put $\overline{\Omega}_+=
\Omega_+\cup\sigma_+\in\Omega(\sigma_1,...,\sigma_{k}, \sigma)$ and
$\overline{\Omega}_-=
\Omega_+\cup\sigma_+\in\Omega(\sigma,\sigma_{k+1},..., \sigma_{m})$.
The simple covering $f:\widetilde{\Omega}\rightarrow\Omega$
generates the simple coverings
$f_{\pm}:(\widetilde{\Omega}_{\pm}\cup f^{-1}(\sigma)_{\pm}))
\rightarrow\overline{\Omega}_{\pm}$, where
$\widetilde{\Omega}_{\pm}= f^{-1}(\Omega_{\pm})$. Consider the set
$M_\Omega^\beta(\beta_1,...,\beta_{m})$ of simple coverings
$f:\widetilde{\Omega}\rightarrow\Omega$, such that local invariants
of $f$ at vertexes of  $\Omega$ are $\beta_i$ and the covering
$f_+|_{f^{-1}(\sigma)_+}$ is equivalent to $\beta$.

The coverings from $M_\Omega^\beta(\beta_1,...,\beta_{m})$ are
coverings that appear from coverings $f_+:(\widetilde{\Omega}_+\cup
f^{-1}(\sigma)_+)) \rightarrow\overline{\Omega}_+$ and
$f_-:(\widetilde{\Omega}_-\cup f^{-1}(\sigma)_-))
\rightarrow\overline{\Omega}_-$ by gluing according to an
equivalence  $f_+|_{f^{-1}(\sigma)_+}\rightarrow
f_-|_{f^{-1}(\sigma)_-}$. Moreover, different equivalences
$f_+|_{f^{-1}(\sigma)_+}\rightarrow f_-|_{f^{-1}(\sigma)_-}$
generate isomorphic coverings if one of it appear from other by
automorphisms of $f_+$ and $f_-$. Thus, $\sum_{f\in
M_\Omega^\beta(\beta_1,...,\beta_{m})}\frac{1}{|\Aut(f)|}=
\sum_{f_+\in\Upsilon_{\overline{\Omega}_+}^d(\beta_1
,...,\beta_{k},\beta), f_-\in
\Upsilon_{\overline{\Omega}_-}^d(\beta,\beta_{k+1},..., \beta_{m})}$
$\frac{|\Aut(\beta)|}{|\Aut(f_+)||\Aut(f_-)|}=$ $\sum_{f_+\in
\Upsilon_{\overline{\Omega}_+}^d(\beta_1 ,...,\beta_{k},\beta)}
\frac{1}{|\Aut(f_+)|}|\Aut(\beta)| \sum_{f_-\in
\Upsilon_{\overline{\Omega}_-}^d(\beta,\beta_{k+1},..., \beta_{m})}
\frac{1}{|\Aut(f_-)|}$. Therefore
$H^d(\Omega,\{\beta_1,...\beta_m\})=
\sum_{\beta\in\Upsilon_{\sigma}^d}\sum_{f\in
M_\Omega^\beta|(\beta_1,...,\beta_{m})} \frac{1}{|\Aut(f)|}=$
$\sum_{\beta\in\Upsilon^d_\beta}
H(\Omega_+,\{\beta_1,...\beta_k,\beta\})|\Aut(\beta)|
H(\Omega_-,\{\beta,\beta_{k+1},...\beta_m\})$. Moreover
$H(\Omega,\{\beta,\gamma\})=\frac{\delta_{\beta,\gamma^*}}{|\Aut(\beta)|}$.

$\Box$

It is follow from theorem \ref{t2.1} that the theory $\mathcal
{H}^d$ generates a Colored Frobenius Algebra $C_*^d$. It defines a
multiplication between isomorphic classes of simple  $d$-coverings
of colored complexes. Theorem \ref{t2.2} gives a formula for
calculation of Hurwitz numbers of simple coverings in terms of these
algebra. The algebra $C_*^d$ is an hight dimensional analog of the
algebra bi-party graphs from \cite{AN1,AN2}. Algebra $C_*^1$ is an
algebra of brane complexes.

\section{Group foam}

\subsection{Hurwitz numbers of $G$-coverings}
Fix a finite group $G$. A homomorphism
$G\rightarrow\Aut(\widetilde{\Omega})$ to the group of automorphism
of simple CW-complex $\widetilde{\Omega}$ is called a
\textit{$G$-action} on $\widetilde{\Omega}$, if the set of fixed
points of any $g\in G$ is sum of some full cell of
$\widetilde{\Omega}$. For brevity we shall identify elements of  $G$
and its images in $\Aut(\widetilde{\Omega})$.

We say that a map $f:\widetilde{\Omega}\rightarrow\Omega$ between
CW-complexes is \textit{G-covering} if $f=f_*r$, where
$r:\widetilde{\Omega}\rightarrow \widetilde{\Omega}/G$ is the
natural projection by some $G$-action and
$f_*:\widetilde{\Omega}/G\rightarrow\Omega$ is an isomorphism.
G-coverings $f':\widetilde{\Omega}'\rightarrow\Omega$ and
$f":\widetilde{\Omega}"\rightarrow\Omega$ are called
\textit{isomorphic}, if there exists an isomorphism
$\widetilde{\varphi}:\widetilde{\Omega}'\rightarrow\widetilde{\Omega}"$
such that $\widetilde{\varphi} g=g\widetilde{\varphi}$, for $g\in
G$, and $f'=f"\widetilde{\varphi}$. Denote by $\Aut^G(f)$ the
automorphisms group of a $G$-covering $f$.

Let $\Omega'$ and $\Omega"$ colored complexes. G-coverings
$f':\widetilde{\Omega}'\rightarrow\Omega'$ and
$f":\widetilde{\Omega}"\rightarrow\Omega"$ are called
\textit{equivalent}, if there exists isomorphisms
$\widetilde{\varphi}:\widetilde{\Omega}'\rightarrow\widetilde{\Omega}"$
and $\varphi:\Omega'\rightarrow\Omega"$ such that
$\widetilde{\varphi} g=g\widetilde{\varphi}$, for $g\in G$, and
$\varphi f'=f"\widetilde{\varphi}$. Denote by $\Upsilon_{\Omega}^G$
the set of equivalent classes of $G$-coverings over a colored
complex $\Omega$.

A $G$-covering $f:\widetilde{\Omega}\rightarrow\Omega$ over a brane
complex $\Omega$ generates a $G$-covering
$f_q:\widetilde{\Omega}_q\rightarrow\sigma_q$ over the vertex
complex $\sigma_q$ of a vertex $q$. The equivalence class of $f_q$
is called a \textit{local invariant} of $f$ at $q$. Consider a set
$\{\beta_q\}=\{\beta_q\in\Upsilon_{\sigma_q}^G|q\in\Omega_b\}$.
Denote by $\Upsilon_{\Omega}^G(\{\beta_q\})$ the set of isomorphic
classes of $G$-coverings $f:\widetilde{\Omega}\rightarrow\Omega$
with local invariants $f_q=\beta_q$ for all vertexes $q\in\Omega_b$.
The weighted number of $G$-covering
$$H^G(\Omega,\{\beta_q\}) = \sum_{f\in
\Upsilon_{\Omega}^G(\{\beta_q\})} 1/|\Aut^G(f)|$$ is called
\textit{Hurwitz Number of $G$-coverings}

\subsection{Brane Topological Field Theory of $G$-coverings}
For any $\sigma\in\Sigma=\Sigma(D,S)$, denote by $C_{\sigma}^G$ the
vector space generated by the set $\Upsilon_{\sigma}^G$. The
involution $\sigma\mapsto\sigma^*$ generates the involution
$C_{\sigma}^G\mapsto C_{\sigma^*}^G$. Hurwitz numbers generate a
family of linear functions $\mathcal {H}^G=
\{\Phi_\Omega:V_{\Omega}\to \mathbb{K}\}$ on
$\{V_{\Omega}=(\otimes_{q\in\Omega_b}C_q^G)\}$, depending from brane
complexes $\Omega$. By repeating the proof of theorem \ref{t3.1} we
obtain

\begin{theorem} \label{t4.1}  $\mathcal{H}^G$ is a Brane Topological Field Theory.
\end{theorem}

It is follow from theorem \ref{t4.1} that $\mathcal {H}^G$ generates
a Colored Frobenius Algebra $C_*^G$ (\textit{$G$-foam}). It is
obviously that its isomorphic class depends only on $G$ and $|S|$.
Thus, $G$-foams are invariants of $G$ that describe the family of
subgroups of $G$. The restriction of $G$-foam on the category of
2-dimensional disks (see remark \ref{note}) were consider in
\cite{AN2}. The algebra $C_*^G$ coincides with the algebra $C_*^d$
from previous section, if $G$ is the symmetrical group $S_d$.

\end{document}